\documentclass[11pt]{amsart}

\usepackage{amsmath}
\usepackage{amsfonts}
\usepackage[active]{srcltx}

\newtheorem{theorem}{Theorem}

\newtheorem{definition}{Definition}
\newtheorem{corollary}{Corollary}

\newtheorem*{Go1}{Theorem G1}
\newtheorem*{Go2}{Theorem G2}
\newtheorem*{Go3}{Theorem G3}
\newtheorem*{Wa1}{Theorem W1}
\newtheorem*{Wa2}{Theorem W2}
\newtheorem*{Sa}{Theorem S}
\newtheorem*{Zh}{Theorem Zh}

\begin{document}
\author{Ushangi Goginava and Artur Sahakian}
\title[convergence of Ces\'aro means ]{On the convergence of Ces\`{a}ro means of negative order of double
trigonometric Fourier series of functions of bounded partial generalized
variation}
\address{U. Goginava, Institute of Mathematics, Faculty of Exact and Natural
Sciences, Tbilisi State University, Chavcha\-vadze str. 1, Tbilisi 0128,
Georgia}
\email{z\_goginava@hotmail.com}
\address{A. Sahakian, Yerevan State University, Faculty of Mathematics and Mechanics,
Alex Manoukian str. 1, Yerevan 0025, Armenia}
\email{sart@ysu.am}
\maketitle

\begin{abstract}
The convergence of Ces\`{a}ro means of negative order of double
trigonometric Fourier series of functions of bounded partial $\Lambda $%
-variation is investigated. The sufficient and neccessary conditions on the
sequence $\Lambda =\{\lambda _{n}\}$ are found for the convergence of
Ces\`{a}ro means of Fourier series of functions of bounded partial $\Lambda $%
-variation.
\end{abstract}

\medskip

\footnotetext{
2000 Mathematics Subject Classification 42B08
}

\section{Classes of Functions of Bounded Generalized Variation}

In 1881 Jordan \cite{Jo} introduced the class of functions of bounded
variation and applied it to the theory of Fourier series. Hereinafter this
notion was generalized by many authors (quadratic variation, $\Phi$%
-variation, $\Lambda$-variation ets., see \cite{Ch}, \cite{M}, \cite{Wa1}, \cite{Wi}). In two
dimensional case the class BV of functions of bounded variation was
introduced by Hardy \cite{Ha}.

Let $f$ be a real function of two variable of period $2\pi $ with respect to
each variable. Given intervals $I=(a,b)$, $J=(c,d)$ and points $x,y$ from $%
T:=[0,2\pi ]$ we denote
\begin{equation*}
f(I,y):=f(b,y)-f(a,y),\qquad f(x,J)=f(x,d)-f(x,c)
\end{equation*}
and
\begin{equation*}
f(I,J):=f(a,c)-f(a,d)-f(b,c)+f(b,d).
\end{equation*}
Let $E=\{I_{i}\}$ be a collection of nonoverlapping intervals from $T$
ordered in arbitrary way and let $\Omega $ be the set of all such
collections $E$.

The Hardy class BV consists of functions $f$ satisfying the condition
$$
\sup_{E\in \Omega }\sum_{i}{|f(I_{i},0)|}+
 \sup_{x}\sup_{F\in \Omega }\sum_{j}{|f(0,J_{j})|}+
\sup_{F,\,E\in \Omega }\sum_{i}\sum_{j}{|f(I_{i},J_{j})|}<\infty,
$$
where $E=\{I_{i}\}$ and $F=\{J_{j}\}$.

In \cite{GoEJA} U. Goginava introduced the class $PBV$ of functions of bounded partial bounded variation, i.e. functions  $f$ having uniformly bounded variation with respect to each variable:
$$
\sup_{y}\sup_{E\in \Omega }\sum_{i}{|f(I_{i},y)|}+
 \sup_{x}\sup_{F\in \Omega }\sum_{j}{|f(x,J_{j})|}<\infty.
$$
For the sequence of positive numbers $\Lambda =\{\lambda
_{n}\}_{n=1}^{\infty }$ we denote
\begin{equation*}
\Lambda V_{1}(f)=\sup_{y}\sup_{E\in \Omega }\sum_{i}\frac{|f(I_{i},y)|}{%
\lambda _{i}}\,\,\,\,\,\,\left( E=\{I_{i}\}\right) ,
\end{equation*}
\begin{equation*}
\Lambda V_{2}(f)=\sup_{x}\sup_{F\in \Omega }\sum_{j}\frac{|f(x,J_{j})|}{%
\lambda _{j}}\qquad (F=\{J_{j}\}),
\end{equation*}
\begin{equation*}
\Lambda V_{1,2}(f)=\sup_{F,\,E\in \Omega }\sum_{i}\sum_{j}\frac{%
|f(I_{i},J_{j})|}{\lambda _{i}\lambda _{j}}.
\end{equation*}

\begin{definition}
We say that the function $f$ has Bounded $\Lambda $-variation on $T^2=[0,2\pi
]^{2}$ and write $f\in \Lambda BV$, if
\begin{equation*}
\Lambda V(f):=\Lambda V_{1}(f)+\Lambda V_{2}(f)+\Lambda V_{1,2}(f)<\infty .
\end{equation*}
We say that the function $f$ has Bounded Partial $\Lambda $-variation and
write $f\in P\Lambda BV$ if
\begin{equation*}
P\Lambda V(f):=\Lambda V_{1}(f)+\Lambda V_{2}(f)<\infty .
\end{equation*}
\end{definition}

If $\lambda_n\equiv 1$ (or if $0<c<\lambda_n<C<\infty,\ n=1,2,\ldots$) the
classes $\Lambda BV$ and $P\Lambda BV$ coincide with the Hardy class $BV$
and PBV respectively. Hence it is reasonable to assume that $%
\lambda_n\to\infty$ and since the intervals in $E=\{I_i\}$ are ordered
arbitrarily, we will suppose, without loss of generality, that the sequence $%
\{\lambda_n\}$ is increasing. Thus,
\begin{equation}  \label{Lambda}
1<\lambda_1\leq \lambda_2\leq\ldots,\qquad \lim_{n\to\infty}\lambda_n=\infty.
\end{equation}

In the case when $\lambda _{n}=n,\ n=1,2\ldots $ we say \textit{Harmonic
Variation} instead of $\Lambda $-variation and write $H$ instead of $\Lambda$
($HBV$, $PHBV$, $HV(f)$, etc).

The notion of $\Lambda $-variation was introduced by D. Waterman \cite{Wa1}
in one dimensional case and A. Sahakian \cite{Saha} in two dimensional case. 

\begin{definition}[Waterman \cite{Wa2}] Let
$\Lambda=\{\lambda_n\}_{n=1}^\infty$ and
$\Lambda_k=\{\lambda_n\}_{n=k}^\infty$, $k=1,2,\ldots$. We say that the function $f$
is continuous in $\Lambda $-variation and
write $f\in C\Lambda BV$, if
\begin{equation*}
\lim_{k\to\infty}\Lambda_k V(f)=0.
\end{equation*}
\end{definition}

\section{$\left( C;\alpha ,\beta \right) \,\left( -1<\alpha ,\beta<0\right) $ Summability}

Let $f\in L^{1}\left( T^{2}\right) .$ The Fourier series of $f$ with respect
to the trigonometric system is the series
\begin{equation*}
S\left[ f,\left( x,y\right) \right] :=\sum_{m,n=-\infty }^{+\infty }\widehat{%
f}\left( m,n\right) e^{imx}e^{iny},
\end{equation*}
where
\begin{equation*}
\widehat{f}\left( m,n\right) =\frac{1}{4\pi ^{2}}\int_{0}^{2\pi
}\int_{0}^{2\pi }f(x,y)e^{-imx}e^{-iny}dxdy
\end{equation*}
are the Fourier coefficients of the function $f$. The rectangular partial
sums are defined as follows:
\begin{equation*}
S_{M,N}\left[ f,(x,y)\right] :=\sum_{m=-M}^{M}\sum_{n=-N}^{N}\widehat{f}%
\left( m,n\right) e^{imx}e^{iny},
\end{equation*}

The Ces\`{a}ro $(C;\alpha ,\beta ),$ $\alpha ,\beta >-1$, means of
two-dimensional Fourier series are defined by
\begin{equation*}
\sigma _{n,m}^{\alpha ,\beta }f(x,y):=\frac{1}{A_{n}^{\alpha }}\frac{1}{%
A_{m}^{\beta }}\sum_{i=0}^{n}\sum_{j=0}^{m}A_{n-i}^{\alpha -1}A_{m-j}^{\beta
-1}S_{i,j}\left[ f,(x,y)\right]
\end{equation*}
where
\begin{equation*}
A_{0}^{\alpha }=1,\,\,\,A_{k}^{\alpha }=\frac{(\alpha +1)\cdots (\alpha +k)}{%
k!},\quad k=1,2,....
\end{equation*}
We say that the double
trigonometric Fourier series of the function $f$ is  $(C;-\alpha
,-\beta )$ summable to $f$, if
$$
\lim_{n,\,m\to \infty}\sigma _{n,m}^{\alpha ,\beta }f(x,y)=f(x,y).
$$
It is well-known that (see \cite{Zy}, p. 157 )
\begin{equation*}
\sigma _{mn}^{\left( \alpha ,\beta \right) }f\left( x,y\right)
=\frac{1}{\pi ^{2}}\int\limits_{-\pi}^{\pi
}\int\limits_{-\pi}^{\pi }f\left( x+t,y+s\right) K_{m}^{\alpha
}\left( s\right) K_{n}^{\beta }\left( t\right) dsdt,
\end{equation*}
where the kernel $K_n^\alpha$, $-1<\alpha<0$ satisfies the following conditions:
\begin{equation}
\left|K_{n}^{-\alpha }( u)\right|\leq2n,\quad u\in T,
\end{equation}
\begin{equation}
K_{n}^{\alpha }\left( u\right) =\varphi _{n}^{\alpha }\left( u\right)
+O\left( 1/nt^{2}\right) ,\,\,\,0\leq |u|\leq \pi,
\end{equation}
where
\begin{equation}
\varphi _{n}^{\alpha }\left( u\right) =\frac{\sin \left[ \left(
n+1/2+\alpha /2\right) u-\alpha \pi /2\right] }{A_{n}^{\alpha }\left[ 2\sin
u/2\right] ^{1+\alpha }},
\end{equation}
The coefficients $A_n^\alpha$ have following bounds:
\begin{equation}
c_{1}(\alpha )n^{\alpha }\leq A_{n}^{\alpha }\leq c_{2}(\alpha )n^{\alpha }.
\end{equation}
Denote
\begin{equation*}
_{1}\Delta _{i}^{m}f\left( x,y\right) :=f\left( x+\frac{2i\pi }{m},\ y\right)
-f\left( x+\frac{\left( 2i+1\right) \pi }{m},y\right) ,
\end{equation*}
\begin{equation*}
_{2}\Delta _{j}^{n}f\left( x,y\right) :=f\left( x,\ y+\frac{2j\pi }{n}\right)
-f\left( x,\ y+\frac{\left( 2j+1\right) \pi }{n}\right) ,
\end{equation*}
\begin{equation*}
\Delta _{ij}^{mn}f\left( x,y\right) =f\left( x+\frac{2i\pi }{m},\ y+\frac{%
2j\pi }{n}\right) -f\left( x+\frac{\left( 2i+1\right) \pi }{m},\ y+\frac{2j\pi
}{n}\right)
\end{equation*}
\begin{equation*}
-f\left( x+\frac{2i\pi }{m},\ y+\frac{\left( 2j+1\right) \pi }{n}\right)
+f\left( x+\frac{\left( 2i+1\right) \pi }{m},\ y+\frac{\left( 2j+1\right) \pi
}{n}\right) .
\end{equation*}

\section{Formulation of Problems}

Let $C(T^{2})$ be the space of $2\pi $%
-periodic with respect to each variable continuous functions with the norm
\begin{equation*}
\Vert f\Vert _{C}:=\sup_{x,y\in T^{2}}|f(x,y)|.
\end{equation*}

For the function $f(x,y)$ we denote by $f\left( x\pm 0,y\pm 0\right) $ the open coordinate quadrant limits (if exist) at the point $%
\left( x,y\right) $ and set
\begin{eqnarray}\label{limits}
&&\qquad\sum f\left( x\pm0,y\pm0\right)\\
&=&\hskip-2mm\big\{f\left( x+0,y+0\right) +f\left( x+0,y-0\right) +f\left(
x-0,y+0\right) +f\left( x-0,y-0\right) \big\}.\notag
\end{eqnarray}
The well known Dirichlet-Jordan theorem (see \cite{Zy}) states that the
Fourier series of a function $f(x),\ x\in T$ of bounded variation converges
at every point $x$ to the value $\left[ f\left( x+0\right) +f\left(
x-0\right) \right] /2$. If $f$ is in addition continuous on $T$ the Fourier
series converges uniformly on $T$. This result was generalized by Waterman
\cite{Wa1}.

\begin{Wa1} [Waterman \cite{Wa1}] If $f\in HBV$, then $S[f]$  converges
at every point $x$ to the value $\left[ f\left( x+0\right) +f\left( x-0\right)
\right] /2$. If $f$ is in addition continuous on $T$, then $S[f]$ converges uniformly on $T$.
\end{Wa1}

Hardy \cite{Ha} generalized the Dirichlet-Jordan theorem to the double
Fourier series. He proved that if function $f(x,y)$ has bounded variation in
the sense of Hardy ($f\in BV$), then $S \left[ f\right] $ converges at any
point $\left( x,y\right) $ to the value $\frac14\sum f\left(
x\pm0,y\pm0\right) $. If $f$ is in addition continuous on $T^{2}$ then $S
\left[ f\right] $ converges uniformly on $T^{2}$.

\begin{Sa} [Sahakian \cite{Saha}] \textit{The Fourier
series of a function $f\left( x,y\right) \in HBV$ converges to $\frac14\sum
f\left( x\pm0,y\pm0\right) $ at any point $\left( x,y\right) $, where the
quadrant limits (\ref{limits}) exist. The convergence is uniformly on any
compact $K$, where the function $f$ is continuous.}
\end{Sa}

Analogs of Theorem S for higher dimensions can be found in \cite{Sab} and
\cite{Bakh}. Convergence of spherical and other partial sums of double
Fourier series of functions of bounded $\Lambda $-variation was investigated
in details by Dyachenko (see \cite{D1}, \cite{D2}, \cite{DW} and references
therein).

The first author \cite{GoEJA} has proved that in Hardy's theorem there is no
need to require the boundedness of mixed variation. In particular, the
following is true

\begin{Go1} [{Goginava} \cite{GoEJA}]
{Let $f\in C\left(
T^{2}\right) \bigcap PBV$. Then $S\left[ f\right] $ converges uniformly on $%
T^{2}$.}
\end{Go1}

For one-dimensional Fourier series Waterman \cite{Wa2} proved the
following

\begin{Wa2}[Waterman \cite{Wa2}] Let $0<\alpha<1$ and $f\in C\{n^{1-\alpha }\}BV$. Then $S[f]$ is everywhere $\left( C,-\alpha \right) $ summable to  the value $\left[ f\left( x+0\right) +f\left( x-0\right)
\right] /2$ and the summability is uniform on each closed interval of continuity.
\end{Wa2}
Later Sablin proved in \cite{Sab}, that for $0<\alpha<1$ the classes $\{n^{1-\alpha }\}BV$ and
$C\{n^{1-\alpha }\}BV$ coincide.

Zhizhiashvili \cite{Zh} has inverstigated the convergence of Ces\`{a}ro means of double trigonometric Fourier series. In particular, the following theorem was proved.

\begin{Zh}[Zhizhiashvili \cite{Zh}] Let $\alpha ,\beta >0$ and $\alpha +\beta <1$. If $f\in BV$, then the double Fourier
series of  $f$ is $\left( C;-\alpha ,-\beta \right) $ summable
to $\frac{1}{4}\sum f\left( x\pm 0,y\pm 0\right) $ in any point $\left(
x,y\right) $ . The convergence is uniformly on any compact $K$, where the
function $f$ is continuous.
\end{Zh}

For functions of partial bounded variation the problem was solved by the first author in  \cite{GoPRI}.

\begin{Go2}
[{Goginava} \cite{GoPRI}] Let $f\in C\left( T^{2}\right) \cap PBV$
 and $\alpha +\beta <1,\ \alpha ,\beta >0.$ Then the double
trigonometric Fourier series of the function $f$ is uniformly $(C;-\alpha
,-\beta )$ summable to $f$.
\end{Go2}

\begin{Go3}
[{Goginava} \cite{GoPRI}]Let $\alpha +\beta \geq 1,\ \alpha ,\beta >0.$ Then
there exists a continuous function $f_{0}\in PBV$ such that the Ces\`{a}ro $%
(C;-\alpha ,-\beta )$ means $\sigma _{n,m}^{-\alpha ,-\beta }\left(
f_{0};0,0\right) \,$of the double trigonometric Fourier series of $f_{0}$
diverge over cubes.
\end{Go3}

In this paper we consider the following problem.
\textit{Let }$\alpha ,\beta \in \left(
0,1\right),\,\alpha +\beta <1.$ \textit{Under what conditions on the sequence }$ \Lambda =\{\lambda _{n}\}$
\textit{the double Fourier series of the function $f\in P\Lambda BV$  is $(C;-\alpha ,-\beta )$ summable.}

The solution is given in Theorems 1 and 2 bellow.

\section{Main Results}

\begin{theorem}
Let $\alpha ,\beta \in \left( 0,1\right) ,\ \alpha
+\beta <1$ and the sequence $\Lambda=\{\lambda_k\}$ satisfies the conditions:
\begin{equation*}
\frac{\lambda _{k}}{k^{1-\left( \alpha +\beta
\right) }}\downarrow 0,\qquad
\sum\limits_{k=1}^{\infty }\frac{\lambda _{k}}{k^{2-\left( \alpha +\beta
\right) }}<\infty.
\end{equation*}
Then the double Fourier series of the function $f\in P\Lambda BV$ is $\left( C;-\alpha
,-\beta \right) $ summable to $\frac{1}{4}\sum f\left( x\pm 0,y\pm 0\right) $
 at any point $\left( x,y\right) $, where the quadrant limits {(\ref {limits})} exist. The convergence is uniform on any compact $K$, where the function $f$ is continuous.
\end{theorem}

\begin{theorem}
Let $\alpha ,\beta \in \left( 0,1\right) ,\,\alpha +\beta <1$ and the sequence $\Lambda=\{\lambda_k\}$ satisfies the conditions:
\begin{equation*}
\frac{\lambda _{k}}{k^{1-\left( \alpha +\beta
\right) }}\downarrow 0,\qquad
\sum\limits_{k=1}^{\infty }\frac{\lambda _{k}}{k^{2-\left( \alpha +\beta
\right) }}=\infty.
\end{equation*}
Then there exists a continuous function $f\in P\Lambda BV$ for which $\left(
C;-\alpha ,-\beta \right) $ means of two-dimensional Fourier series diverges
over cubes at $\left( 0,0\right).$
\end{theorem}

\begin{corollary}
 Let $\alpha ,\beta \in \left( 0,1\right) ,\,\alpha +\beta <1$. \\
a)If $f\in P\left\{ \frac{n^{1-\left( \alpha +\beta \right) }%
}{\log ^{1+\varepsilon }\left( n+1\right) }\right\} BV$ for some $%
\varepsilon >0$, then the double Fourier series of the function $f$ is $%
\left( C;-\alpha ,-\beta \right) $ summable to $\frac{1}{4}\sum f\left( x\pm
0,y\pm 0\right) $ in any point $\left( x,y\right) $ , where the quadrant limits {(\ref{limits})} exist. The convergence is uniform on any
compact $K$, where the function $f$ is continuous.

b) There exists a continuous function $f\in P\left\{ \frac{n^{1-\left( \alpha
+\beta \right) }}{\log \left( n+1\right) }\right\} BV$ such that $\left(
C;-\alpha ,-\beta \right) $ means of two-dimensional Fourier series of $f$ diverges
over cubes at $\left( 0,0\right) .$
\end{corollary}

\begin{corollary}
Let $\alpha ,\beta \in \left( 0,1\right) ,\,\alpha +\beta <1$\thinspace $\,$%
\thinspace and$\,\,f\in PBV$ . Then the double Fourier series of the
function $f$ is $\left( C;-\alpha ,-\beta \right) $ summable to $\frac{1}{4}%
\sum f\left( x\pm 0,y\pm 0\right) $ in any point $\left( x,y\right) $ ,
where the quadratic limits {(\ref{limits})} exist. The convergence is
uniform on any compact $K$, where the function $f$ is continuous.
\end{corollary}

\section{Proofs}

\begin{proof}[Proof of Theorem 1]
It is easy to show that
\begin{equation*}
\sigma _{mn}^{\left( -\alpha ,-\beta \right) }f\left( x,y\right) -\frac{1}{4}%
\sum f\left( x\pm 0,y\pm 0\right)
\end{equation*}
\begin{equation*}
=\frac{1}{\pi ^{2}}\sum\limits_{i=1}^{4}\int\limits_{0}^{\pi
}\int\limits_{0}^{\pi }\varphi _{i}\left( x,y,s,t\right) K_{m}^{-\alpha
}\left( s\right) K_{n}^{-\beta }\left( t\right) dsdt
\end{equation*}
\begin{equation*}
=:\sum\limits_{i=1}^{4}I_{mn}^{\left( k\right) }\left( x,y\right) .
\end{equation*}
where
\begin{equation*}
\varphi _{1}\left( x,y,s,t\right) :=f\left( x+s,y+t\right) -f\left(
x+0,y+0\right),
\end{equation*}
\begin{equation*}
\varphi _{2}\left( x,y,s,t\right) :=f\left( x-s,y+t\right) -f\left(
x-0,y+0\right),
\end{equation*}
\begin{equation*}
\varphi _{3}\left( x,y,s,t\right) :=f\left( x+s,y-t\right) -f\left(
x+0),y-0\right),
\end{equation*}
\begin{equation*}
\varphi _{4}\left( x,y,s,t\right) :=f\left( x-s,y-t\right) -f\left(
x-0,y-0\right).
\end{equation*}

For $I_{mn}^{\left( 1\right) }\left( x,y\right) $ we can write
\begin{equation}
\pi ^{2}I_{mn}^{\left( 1\right) }\left( x,y\right)
\end{equation}
\begin{equation*}
=\left( \int\limits_{0}^{\pi /m}\int\limits_{0}^{\pi
/n}+\int\limits_{0}^{\pi /m}\int\limits_{\pi /n}^{\pi }+\int\limits_{\pi
/m}^{\pi }\int\limits_{0}^{\pi /n}+\int\limits_{\pi /m}^{\pi
}\int\limits_{\pi /n}^{\pi }\right) \left( \varphi _{1}\left( x,y,s,t\right)
K_{m}^{-\alpha }\left( s\right) K_{n}^{-\beta }\left( t\right) dsdt\right)
\end{equation*}
\begin{equation*}
=:\sum\limits_{k=1}^{4}I_{mn}^{\left( 1k\right) }\left( x,y\right) .
\end{equation*}

From (2) we have
\begin{equation}
\left| I_{mn}^{\left( 11\right) }\left( x,y\right) \right| \leq c\left(
\alpha ,\beta \right) mn\int\limits_{0}^{\pi /m}\int\limits_{0}^{\pi
/n}\left| \varphi _{1}\left( x,y,s,t\right) \right| dsdt
\end{equation}
\begin{equation*}
\leq c\left( \alpha ,\beta \right) \sup\limits_{0<s<\pi /m,0<t<\pi /n}\left|
\varphi _{1}\left( x,y,s,t\right) \right| =o\left( 1\right) \,\,\,\,\text{%
as\thinspace \thinspace \thinspace }m,n\rightarrow \infty .
\end{equation*}

Using (3), we obtain
\begin{equation}
I_{mn}^{\left( 12\right) }\left( x,y\right) =\int\limits_{0}^{\pi
/m}\int\limits_{\pi /n}^{\pi }\varphi _{1}\left( x,y,s,t\right)
K_{m}^{-\alpha }\left( s\right) \varphi _{n}^{-\beta }\left( t\right) dsdt
\end{equation}
\begin{equation*}
+\int\limits_{0}^{\pi /m}\int\limits_{\pi /n}^{\pi }\varphi _{1}\left(
x,y,s,t\right) K_{m}^{-\alpha }\left( s\right) O\left( \frac{1}{nt^{2}}%
\right) dsdt
\end{equation*}
\begin{equation*}
=:I_{mn}^{\left( 121\right) }\left( x,y\right) +I_{mn}^{\left( 122\right)
}\left( x,y\right) .
\end{equation*}
We can write
\begin{eqnarray}
&&\left| I_{mn}^{\left( 122\right) }\left( x,y\right) \right| \\
&\leq &\int\limits_{0}^{\pi /m}\int\limits_{\pi /n}^{\pi /\sqrt{n}}\left|
\varphi _{1}\left( x,y,s,t\right) \right| |K_{m}^{-\alpha }\left( s\right)
|O\left( \frac{1}{nt^{2}}\right) dsdt  \notag \\
&&+\int\limits_{0}^{\pi /m}\int\limits_{\pi /\sqrt{n}}^{\pi }\left| \varphi
_{1}\left( x,y,s,t\right) \right| |K_{m}^{-\alpha }\left( s\right) |O\left(
\frac{1}{nt^{2}}\right) dsdt  \notag \\
&\leq &c\left( \alpha ,\beta ,f\right) \left\{ \sup\limits_{0<s<\pi
/m,0<t<\pi /\sqrt{n}}\left| \varphi _{1}\left( x,y,s,t\right) \right|
+\int\limits_{\pi /\sqrt{n}}^{\pi }\frac{dt}{nt^{2}}\right\}  \notag \\
&=&o\left( 1\right) \,\,\,\text{as\thinspace \thinspace \thinspace }%
n,m\rightarrow \infty .  \notag
\end{eqnarray}

In order to estimate $I_{mn}^{\left( 121\right) }\left( x,y\right) $ it is
enough to estimate the following expression
\begin{equation*}
J_{mn}\left( x,y\right) :=n^{\beta }\int\limits_{0}^{\pi /m}\int\limits_{\pi
/n}^{\pi }\varphi _{1}\left( x,y,s,t\right) K_{m}^{-\alpha }\left( s\right)
w_{\beta }\left( t\right) \sin ntdsdt,
\end{equation*}
where
\begin{equation*}
w_{\beta }\left( t\right) =\frac{\cos \frac{1-\beta }{2}t}{\left( \sin
t/2\right) ^{1-\beta }}.
\end{equation*}

We have
\begin{equation*}
J_{mn}\left( x,y\right) =n^{\beta
}\sum\limits_{i=1}^{n-1}\int\limits_{0}^{\pi /m}K_{m}^{-\alpha }\left(
s\right) \left( \int\limits_{i\pi /n}^{\left( i+1\right) \pi /n}\varphi
_{1}\left( x,y,s,t\right) w_{\beta }\left( t\right) \sin ntdt\right) ds
\end{equation*}
\begin{equation*}
=n^{\beta }\sum\limits_{i=1}^{\left( n-1\right) /2}\int\limits_{0}^{\pi
/m}K_{m}^{-\alpha }\left( s\right) \left( \int\limits_{0}^{\pi /n}\left[
\varphi _{1}\left( x,y,s,t+\frac{2i\pi }{n}\right) -\varphi _{1}\left(
x,y,s,t+\frac{\left( 2i+1\right)\pi }{n}\right) \right]\right.
\end{equation*}
\begin{equation*}
\left. w_{\beta }\left( t+\frac{2i\pi }{n}\right) \sin ntdt\right) ds
\end{equation*}
\begin{equation*}
+n^{\beta }\sum\limits_{i=1}^{\left( n-1\right) /2}\int\limits_{0}^{\pi
/m}K_{m}^{-\alpha }\left( s\right) \left( \int\limits_{0}^{\pi /n}\varphi
_{1}\left( x,y,s,t+\frac{\left( 2i+1\right) \pi }{n}\right) \right.
\end{equation*}
\begin{equation*}
\left. \left[ w_{\beta }\left( t+\frac{2i\pi }{n}\right) -w_{\beta }\left( t+%
\frac{\left( 2i+1\right) \pi }{n}\right) \right] \sin ntdt\right) ds
\end{equation*}
\begin{equation*}
=:J_{mn}^{\left( 1\right) }\left( x,y\right) +J_{mn}^{\left( 2\right) }\left(
x,y\right) .
\end{equation*}
Using the following inequality:
\begin{equation}
\left| w_{\beta }\left( t+\frac{2i\pi }{n}\right) -w_{\beta }\left( t+\frac{%
\left( 2i+1\right) \pi }{n}\right) \right| \leq \frac{c\left( \beta \right)
n^{1-\beta }}{i^{2-\beta }},
\end{equation}
for $J_{mn}^{\left( 2\right) }\left( x,y\right) $ we can write
\begin{equation}
\left| J_{mn}^{\left( 2\right) }\left( x,y\right) \right|
\end{equation}
\begin{equation*}
\leq c\left( \beta \right) mn\sum\limits_{i=1}^{\left( n-1\right) /2}\frac{1%
}{i^{2-\beta }}\int\limits_{0}^{\pi /m}\left( \int\limits_{0}^{\pi /n}\left|
\varphi _{1}\left( x,y,s,t+\frac{\left( 2i+1\right) \pi }{n}\right) \right|
dt\right) ds
\end{equation*}
\begin{equation*}
\leq c\left( \beta \right) nm\sum\limits_{i\leq \sqrt{n}}\frac{1}{i^{2-\beta
}}\int\limits_{0}^{\pi /m}\left( \int\limits_{0}^{\pi /n}\left| \varphi
_{1}\left( x,y,s,t+\frac{\left( 2i+1\right) \pi }{n}\right) \right|
dt\right) ds
\end{equation*}
\begin{equation*}
+c\left( \beta \right) nm\sum\limits_{\sqrt{n}<i\leq \left( n-1\right) /2}%
\frac{1}{i^{2-\beta }}\int\limits_{0}^{\pi /m}\left( \int\limits_{0}^{\pi
/n}\left| \varphi _{1}\left( x,y,s,t+\frac{\left( 2i+1\right) \pi }{n}%
\right) \right| dt\right) ds
\end{equation*}
\begin{equation*}
\leq c\left( \beta \right) \sup\limits_{0<s<\pi /n,0<s<4\pi /\sqrt{n}}\left|
\varphi _{1}\left( x,y,s,t\right) \right| +c\left( \beta ,f\right) \left(
\frac{1}{\sqrt{n}}\right) ^{1-\beta }=o\left( 1\right),
\end{equation*}
as $n,m\rightarrow \infty$.

To estimate $J_{mn}^{\left( 1\right) }\left( x,y\right) $, we denote
\begin{equation}
\mu \left( n,m\right) :=\left[\min \left\{ \frac{1}{2}\ln n-1,\ \left( s\left(
n,m\right) \right) ^{-1}\right\}\right],
\end{equation}
where $[a]$ is the integer part of $a$ and
\begin{equation}
s\left( n,m\right) :=\sup\limits_{0<s<\pi /m,\ 0<t<\pi \ln n/n}\left| \varphi
_{1}\left( x,y,s,t\right) \right| .
\end{equation}
Then we have
\begin{equation}
\left| J_{mn}^{\left( 1\right) }\left( x,y\right) \right| \leq c\left( \beta
\right) nm\int\limits_{0}^{\pi /m}\left( \int\limits_{0}^{\pi
/n}\sum\limits_{i=1}^{\mu \left( n,m\right) }\frac{1}{i^{1-\beta }}\left|
\varphi _{1}\left( x,y,s,t+\frac{2i\pi }{n}\right) \right. \right.
\end{equation}
\begin{equation*}
\left. \left. -\varphi _{1}\left( x,y,s,t+\frac{\left( 2i+1\right) \pi }{n}%
\right) \right| dtds\right)
\end{equation*}
\begin{equation*}
+c\left( \beta \right) nm\int\limits_{0}^{\pi /m}\left( \int\limits_{0}^{\pi
/n}\sum\limits_{i=\mu \left( n,m\right) }^{\left( n-1\right) /2}\frac{1}{%
i^{1-\beta }}\left| \varphi _{1}\left( x,y,s,t+\frac{2i\pi }{n}\right)
\right. \right.
\end{equation*}
\begin{equation*}
\left. \left. -\varphi _{1}\left( x,y,s,t+\frac{\left( 2i+1\right) \pi }{n}%
\right) \right| dt\right) ds
\end{equation*}
\begin{equation*}
\leq c\left( \beta \right) \sup\limits_{0<s<\pi /m,\ 0<t<\left( 2\mu \left(n,\,m\right) +1\right) \pi /n}\left| \varphi _{1}\left( x,y,s,t\right) \right|
\left( \mu \left( n,m\right) \right) ^{\beta }
\end{equation*}
\begin{equation*}
+c\left( \beta \right) nm\int\limits_{0}^{\pi /m}\left( \int\limits_{0}^{\pi
/n}\sum\limits_{i=\mu \left( n,m\right) }^{\left( n-1\right) /2}\frac{1}{%
\lambda _{i}}\left| \varphi _{1}\left( x,y,s,t+\frac{2i\pi }{n}\right)
\right. \right.
\end{equation*}
\begin{equation*}
\left. \left. -\varphi _{1}\left( x,y,s,t+\frac{\left( 2i+1\right) \pi }{n}%
\right) \right| \frac{\lambda _{i}}{i^{1-\beta }}dt\right) ds
\end{equation*}
\begin{equation*}
\leq c\left( \beta \right) \sup\limits_{0<s<\pi /m,\ 0<t<\pi \ln n/n}\left|
\varphi _{1}\left( x,y,s,t\right) \right| \left( \mu \left( n,m\right)
\right) ^{\beta }
\end{equation*}
\begin{equation*}
+c\left( \beta \right) nm\frac{\lambda _{\mu \left( n,m\right) }}{\left( \mu
\left( n,m\right) \right) ^{1-\beta }}\int\limits_{0}^{\pi /m}\left(
\int\limits_{0}^{\pi /n}\sum\limits_{i=\mu \left( n,m\right) }^{\left(
n-1\right) /2}\frac{1}{\lambda _{i}}\left| \varphi _{1}\left( x,y,s,t+\frac{%
2i\pi }{n}\right) \right. \right.
\end{equation*}
\begin{equation*}
\left. \left. -\varphi _{1}\left( x,y,s,t+\frac{\left( 2i+1\right) \pi }{n}%
\right) \right| dt\right) ds
\end{equation*}
\begin{equation*}
\leq c\left( \beta \right) s\left( n,m\right) \left( \mu \left( n,m\right)
\right) ^{\beta }+c\left( \beta \right) \frac{\lambda _{\mu \left(
n,m\right) }}{\left( \mu \left( n,m\right) \right) ^{1-\beta }}V_{2}\Lambda
\left( f\right) =o\left( 1\right) \,\,\,\,\text{as\thinspace \thinspace
\thinspace }n,m\rightarrow \infty .
\end{equation*}

Combining (9), (10), (12) and (15) we conclude that
\begin{equation}
I_{mn}^{\left( 12\right) }\left( x,y\right) \rightarrow 0\,\,\,\text{%
as\thinspace \thinspace \thinspace }m,n\rightarrow \infty .
\end{equation}

Analogously, we can prove that
\begin{equation}
I_{mn}^{\left( 13\right) }\left( x,y\right) \rightarrow 0\,\,\,\text{%
as\thinspace \thinspace \thinspace }m,n\rightarrow \infty .
\end{equation}

In order to estimate $I_{mn}^{\left( 14\right) }\left( x,y\right) $ it is
enough to estimate the following expression
\begin{equation*}
L_{mn}\left( x,y\right) :=m^{\alpha }n^{\beta }\int\limits_{\pi /m}^{\pi
}\int\limits_{\pi /n}^{\pi }\varphi _{1}\left( x,y,s,t\right) w_{\alpha
}\left( s\right) w_{\beta }\left( t\right) \sin ms\sin ntdtds.
\end{equation*}
We have
\begin{equation}
L_{mn}\left( x,y\right) =m^{\alpha }n^{\beta }\sum\limits_{i=1}^{\left(
m-1\right) /2}\sum\limits_{j=1}^{\left( n-1\right) /2}\int\limits_{0}^{\pi
/m}\int\limits_{0}^{\pi /n}\varphi _{1}\left( x,y,s+\frac{2i\pi }{m},t+\frac{%
2j\pi }{n}\right)
\end{equation}
\begin{equation*}
\times w_{\alpha }\left( s+\frac{2i\pi }{m}\right) w_{\beta }\left( t+\frac{%
2j\pi }{n}\right) \sin ms\sin ntdtds
\end{equation*}
\begin{equation*}
-m^{\alpha }n^{\beta }\sum\limits_{i=1}^{\left( m-1\right)
/2}\sum\limits_{j=1}^{\left( n-1\right) /2}\int\limits_{0}^{\pi
/m}\int\limits_{0}^{\pi /n}\varphi _{1}\left( x,y,s+\frac{\left( 2i+1\right)
\pi }{m},t+\frac{2j\pi }{n}\right)
\end{equation*}
\begin{equation*}
\times w_{\alpha }\left( s+\frac{\left( 2i+1\right) \pi }{m}\right) w_{\beta
}\left( t+\frac{2j\pi }{n}\right) \sin ms\sin ntdtds
\end{equation*}
\begin{equation*}
-m^{\alpha }n^{\beta }\sum\limits_{i=1}^{\left( m-1\right)
/2}\sum\limits_{j=1}^{\left( n-1\right) /2}\int\limits_{0}^{\pi
/m}\int\limits_{0}^{\pi /n}\varphi _{1}\left( x,y,s+\frac{2i\pi }{m},t+\frac{%
\left( 2j+1\right) \pi }{n}\right)
\end{equation*}
\begin{equation*}
\times w_{\alpha }\left( s+\frac{2i\pi }{m}\right) w_{\beta }\left( t+\frac{%
\left( 2j+1\right) \pi }{n}\right) \sin ms\sin ntdtds
\end{equation*}
\begin{equation*}
+m^{\alpha }n^{\beta }\sum\limits_{i=1}^{\left( m-1\right)
/2}\sum\limits_{j=1}^{\left( n-1\right) /2}\int\limits_{0}^{\pi
/m}\int\limits_{0}^{\pi /n}\varphi _{1}\left( x,y,s+\frac{\left( 2i+1\right)
\pi }{m},t+\frac{\left( 2j+1\right) \pi }{n}\right)
\end{equation*}
\begin{equation*}
\times w_{\alpha }\left( s+\frac{\left( 2i+1\right) \pi }{m}\right) w_{\beta
}\left( t+\frac{\left( 2j+1\right)\pi }{n}\right) \sin ms\sin ntdtds
\end{equation*}
\begin{equation*}
=m^{\alpha }n^{\beta }\sum\limits_{i=1}^{\left( m-1\right)
/2}\sum\limits_{j=1}^{\left( n-1\right) /2}\int\limits_{0}^{\pi
/m}\int\limits_{0}^{\pi /n}\left[ \varphi _{1}\left( x,y,s+\frac{2i\pi }{m}%
,t+\frac{2j\pi }{n}\right) \right.
\end{equation*}
\begin{equation*}
-\varphi _{1}\left( x,y,s+\frac{\left( 2i+1\right) \pi }{m},t+\frac{2j\pi }{n%
}\right) -\varphi _{1}\left( x,y,s+\frac{2i\pi }{m},t+\frac{\left(
2j+1\right) \pi }{n}\right)
\end{equation*}
\begin{equation*}
\left. +\varphi _{1}\left( x,y,s+\frac{\left( 2i+1\right) \pi }{m},t+\frac{%
\left( 2j+1\right) \pi }{n}\right) \right]
\end{equation*}
\begin{equation*}
\times w_{\alpha }\left( s+\frac{2i\pi }{m}\right) w_{\beta }\left( t+\frac{%
2j\pi }{n}\right) \sin ms\sin ntdtds
\end{equation*}

\begin{equation*}
+m^{\alpha }n^{\beta }\sum\limits_{i=1}^{\left( m-1\right)
/2}\sum\limits_{j=1}^{\left( n-1\right) /2}\int\limits_{0}^{\pi
/m}\int\limits_{0}^{\pi /n}\left[ \varphi _{1}\left( x,y,s+\frac{\left(
2i+1\right) \pi }{m},t+\frac{2j\pi }{n}\right) \right.
\end{equation*}
\begin{equation*}
\left. -\varphi _{1}\left( x,y,s+\frac{\left( 2i+1\right) \pi }{m},t+\frac{%
\left( 2j+1\right) \pi }{n}\right) \right]
\end{equation*}
\begin{equation*}
\times \left[ w_{\alpha }\left( s+\frac{2i\pi }{m}\right) -w_{\alpha }\left(
s+\frac{\left( 2i+1\right) \pi }{m}\right) \right] w_{\beta }\left( t+\frac{%
2j\pi }{n}\right) \sin ms\sin ntdtds
\end{equation*}
\begin{equation*}
+m^{\alpha }n^{\beta }\sum\limits_{i=1}^{\left( m-1\right)
/2}\sum\limits_{j=1}^{\left( n-1\right) /2}\int\limits_{0}^{\pi
/m}\int\limits_{0}^{\pi /n}\left[ \varphi _{1}\left( x,y,s+\frac{2i\pi }{m}%
,t+\frac{\left( 2j+1\right) \pi }{n}\right) \right.
\end{equation*}
\begin{equation*}
\left. -\varphi _{1}\left( x,y,s+\frac{\left( 2i+1\right) \pi }{m},t+\frac{%
\left( 2j+\right) 1\pi }{n}\right) \right]
\end{equation*}
\begin{equation*}
\times \left[ w_{\beta }\left( t+\frac{2j\pi }{n}\right) -w_{\beta }\left( t+%
\frac{\left( 2j+1\right) \pi }{n}\right) \right] w_{\alpha }\left( s+\frac{%
2i\pi }{m}\right) \sin ms\sin ntdtds
\end{equation*}
\begin{equation*}
+m^{\alpha }n^{\beta }\sum\limits_{i=1}^{\left( m-1\right)
/2}\sum\limits_{j=1}^{\left( n-1\right) /2}\int\limits_{0}^{\pi
/m}\int\limits_{0}^{\pi /n}\varphi _{1}\left( x,y,s+\frac{\left( 2i+1\right)
\pi }{m},t+\frac{\left( 2j+1\right) \pi }{n}\right)
\end{equation*}
\begin{equation*}
\times \left[ w_{\beta }\left( t+\frac{2j\pi }{n}\right) -w_{\beta }\left( t+%
\frac{\left( 2j+1\right) \pi }{n}\right) \right]
\end{equation*}
\begin{equation*}
\left[ w_{\alpha }\left( s+\frac{2i\pi }{m}\right) -w_{\alpha }\left( s+%
\frac{\left( 2i+1\right) \pi }{m}\right) \right] \sin ms\sin ntdtds
\end{equation*}
\begin{equation*}
=:\sum\limits_{k=1}^{4}L_{mn}^{\left( k\right) }\left( x,y\right) .
\end{equation*}

By (11) we obtain
\begin{equation}
\left| L_{mn}^{\left( 4\right) }\left( x,y\right) \right| \leq c\left(
\alpha ,\beta \right) mn\sum\limits_{i=1}^{\left[ \sqrt{m}\right] }\frac{1}{%
i^{2-\alpha }}\sum\limits_{j=1}^{\left[ \sqrt{n}\right] }\frac{1}{j^{2-\beta
}}
\end{equation}
\begin{equation*}
\int\limits_{0}^{\pi /m}\int\limits_{0}^{\pi /n}\left| \varphi _{1}\left(
x,y,s+\frac{\left( 2i+1\right) \pi }{m},t+\frac{\left( 2j+1\right) \pi }{n}%
\right) \right|
\end{equation*}
\begin{equation*}
+c\left( \alpha ,\beta ,f\right) \sum\limits_{i=1}^{\infty }\frac{1}{%
i^{2-\alpha }}\sum\limits_{j=\left[ \sqrt{n}\right] }^{\infty }\frac{1}{%
j^{2-\beta }}
\end{equation*}
\begin{equation*}
+c\left( \alpha ,\beta ,f\right) mn\sum\limits_{i=\left[ \sqrt{m}\right]
}^{\infty }\frac{1}{i^{2-\alpha }}\sum\limits_{j=1}^{\infty }\frac{1}{%
j^{2-\beta }}
\end{equation*}
\begin{equation*}
\leq c\left( \alpha ,\beta \right) \sup\limits_{0<s<4\pi /\sqrt{m},0<t<4\pi /%
\sqrt{n}}\left| \varphi _{1}\left( x,y,s,t\right) \right| +o\left( 1\right)
\end{equation*}
\begin{equation*}
=o\left( 1\right) \text{\thinspace \thinspace \thinspace \thinspace as }%
n,m\rightarrow \infty .
\end{equation*}

Let
\begin{equation}
\tau \left( n,m\right) :=\left[\min \left\{ \frac{1}{2}\ln n-1,\frac{1}{2}\ln
m-1,\left( l\left( n,m\right) \right) ^{-1}\right\}\right],
\end{equation}
where
\begin{equation*}
l\left( n,m\right) :=\sup\limits_{0<s<\pi \ln m/m,\ 0<t<\pi \ln n/n}\left|
\varphi _{1}\left( x,y,s,t\right) \right|
\end{equation*}

Then we can write
\begin{equation}
\left| L_{mn}^{\left( 3\right) }\left( x,y\right) \right|
\end{equation}
\begin{equation*}
\leq c\left( \alpha ,\beta \right) mn\sum\limits_{i=1}^{\tau \left(
n,m\right) }\frac{1}{i^{1-\alpha }}\sum\limits_{j=1}^{\tau \left( n,m\right)
}\frac{1}{j^{2-\beta }}\int\limits_{0}^{\pi /m}\int\limits_{0}^{\pi
/n}\left| \varphi _{1}\left( x,y,s+\frac{2i\pi }{m},t+\frac{\left(
2j+1\right) \pi }{n}\right) \right.
\end{equation*}
\begin{equation*}
\left. -\varphi _{1}\left( x,y,s+\frac{\left( 2i+1\right) \pi }{m},t+\frac{%
\left( 2j+1\right) \pi }{n}\right) \right| dtds
\end{equation*}
\begin{equation*}
+c\left( \alpha ,\beta \right) mn\sum\limits_{i=\tau \left( n,m\right)
}^{\left( m-1\right) /2}\frac{1}{\lambda _{i}}\frac{\lambda _{i}}{%
i^{1-\alpha }}\sum\limits_{j=1}^{\tau \left( n,m\right) }\frac{1}{j^{2-\beta
}}\int\limits_{0}^{\pi /m}\int\limits_{0}^{\pi /n}\left| \varphi _{1}\left(
x,y,s+\frac{2i\pi }{m},t+\frac{\left( 2j+1\right) \pi }{n}\right) \right.
\end{equation*}
\begin{equation*}
\left. -\varphi _{1}\left( x,y,s+\frac{\left( 2i+1\right) \pi }{m},t+\frac{%
\left( 2j+1\right) \pi }{n}\right) \right| dtds
\end{equation*}
\begin{equation*}
+c\left( \alpha ,\beta \right) mn\sum\limits_{i=1}^{\left( m-1\right) /2}%
\frac{1}{\lambda _{i}}\frac{\lambda _{i}}{i^{1-\alpha }}\sum\limits_{j=\tau
\left( n,m\right) }^{\left( n-1\right) /2}\frac{1}{j^{2-\beta }}%
\int\limits_{0}^{\pi /m}\int\limits_{0}^{\pi /n}\left| \varphi _{1}\left(
x,y,s+\frac{2i\pi }{m},t+\frac{\left( 2j+1\right) \pi }{n}\right) \right.
\end{equation*}
\begin{equation*}
\left. -\varphi _{1}\left( x,y,s+\frac{\left( 2i+1\right) \pi }{m},t+\frac{%
\left( 2j+1\right) \pi }{n}\right) \right| dtds
\end{equation*}
\begin{equation*}
\leq c\left( \alpha ,\beta \right) l\left( n,m\right) \left( \tau \left(
n,m\right) \right) ^{\alpha +\beta }
\end{equation*}
\begin{equation*}
+c\left( \alpha ,\beta \right) \frac{\lambda _{\tau \left( n,m\right) }}{%
\left( \tau \left( n,m\right) \right) ^{1-\alpha }}V_{1}\Lambda \left(
f\right) +c\left( \alpha ,\beta \right) \frac{1}{\left( \tau \left(
n,m\right) \right) ^{1-\beta }}V_{1}\Lambda \left( f\right)
\end{equation*}
\begin{equation*}
=o\left( 1\right) \,\,\,\text{as\thinspace \thinspace \thinspace }%
n,m\rightarrow \infty .
\end{equation*}

Analogously, we can prove that
\begin{equation}
\left| L_{mn}^{\left( 2\right) }\left( x,y\right) \right| =o\left( 1\right)
\,\,\,\text{as\thinspace \thinspace \thinspace }n,m\rightarrow \infty .
\end{equation}

For $L_{mn}^{\left( 1\right) }\left( x,y\right) $ we can write
\begin{equation}
\left| L_{mn}^{\left( 1\right) }\left( x,y\right) \right|
\end{equation}
\begin{equation*}
\leq c\left( \alpha ,\beta \right) mn\left\{ \sum\limits_{i=1}^{\tau \left(
n,m\right) }\frac{1}{i^{1-\alpha }}\sum\limits_{j=1}^{\tau \left( n,m\right)
}\frac{1}{j^{1-\beta }}+\sum\limits_{i=\tau \left( n,m\right) }^{\left(
m-1\right) /2}\frac{1}{i^{1-\alpha }}\sum\limits_{j=1}^{\tau \left(
n,m\right) }\frac{1}{j^{1-\beta }}\right.
\end{equation*}
\begin{equation*}
\left. +\sum\limits_{i=1}^{\tau \left( n,m\right) }\frac{1}{i^{1-\alpha }}%
\sum\limits_{j=\tau \left( n,m\right) }^{\left( n-1\right) /2}\frac{1}{%
j^{1-\beta }}+\sum\limits_{i=\tau \left( n,m\right) }^{\left( m-1\right) /2}%
\frac{1}{i^{1-\alpha }}\sum\limits_{j=\tau \left( n,m\right) }^{\left(
n-1\right) /2}\frac{1}{j^{1-\beta }}\right\}
\end{equation*}
\begin{equation*}
\left( \int\limits_{0}^{\pi /m}\int\limits_{0}^{\pi /n}|\Delta
_{ij}^{mn}f\left( x+s,y+t\right) |dsdt\right)
=:\sum\limits_{k=1}^{4}L_{mn}^{\left( 1k\right) }\left( x,y\right) .
\end{equation*}

From (20) we obtain that
\begin{equation}
\left| L_{mn}^{\left( 11\right) }\left( x,y\right) \right| \leq c\left(
\alpha ,\beta \right) \left( l\left( n,m\right) ^{1-\alpha -\beta }\right)
\rightarrow 0\,\,\,\text{as\thinspace \thinspace \thinspace }m,n\rightarrow
\infty .
\end{equation}
Next, we have
\begin{equation}
\left| L_{mn}^{\left( 13\right) }\left( x,y\right) \right| \leq c\left(
\alpha ,\beta \right) mn\left\{ \sum\limits_{i=1}^{\tau \left( n,m\right) }%
\frac{1}{i^{1-\alpha }}\sum\limits_{j=\tau \left( n,m\right) }^{\left(
n-1\right) /2}\frac{1}{j^{1-\beta }}\right.
\end{equation}
\begin{equation*}
\left. \left( \int\limits_{0}^{\pi /m}\int\limits_{0}^{\pi /n}|\Delta
_{ij}^{mn}f\left( x+s,y+t\right) |dsdt\right) \right\}
\end{equation*}
\begin{equation*}
\leq c\left( \alpha ,\beta \right) n\left\{ \sum\limits_{i=1}^{\tau \left(
n,m\right) }\frac{1}{i^{1-\alpha }}\left( \int\limits_{0}^{\pi /n}\sup\limits_{x}\sum\limits_{j=\tau \left(
n,m\right) }^{\left( n-1\right) /2}\frac{\lambda _{j}}{j^{1-\beta }}\frac{1}{%
\lambda _{j}}|_{2}\Delta _{j}^{n}f\left(
x,y+t\right) |dt\right) \right\}
\end{equation*}
\begin{equation*}
\leq c\left( \alpha ,\beta \right) \frac{\tau \left( n,m\right) }{\left(
\tau \left( n,m\right) \right) ^{1-\beta -\alpha }}V_{2}\Lambda \left(
f\right) \rightarrow 0,
\quad\text{as}\ n,m\rightarrow \infty .
\end{equation*}

Analogously, we can prove that
\begin{equation}
\left| L_{mn}^{\left( 12\right) }\left( x,y\right) \right| \rightarrow 0,
\quad\text{as}\ n,m\rightarrow \infty .
\end{equation}

From the condition of the Theorem 1 we can write
\begin{equation}
\left| L_{mn}^{\left( 14\right) }\left( x,y\right) \right|
\end{equation}
\begin{equation*}
\leq c\left( \alpha ,\beta \right) nm\sum\limits_{i=\tau \left( n,m\right)
}^{\left( m-1\right) /2}\frac{1}{i^{1-\alpha }}\sum\limits_{j=\tau \left(
n,m\right) }^{\left( n-1\right) /2}\frac{1}{j^{1-\beta }}
\end{equation*}
\begin{equation*}
\left( \int\limits_{0}^{\pi /m}\int\limits_{0}^{\pi /n}|\Delta
_{ij}^{mn}f\left( x+s,y+t\right) |dsdt\right)
\end{equation*}
\begin{equation*}
\leq c\left( \alpha ,\beta \right) nm\left\{ \sum\limits_{i=\tau \left(
n,m\right) }^{\left( m-1\right) /2}\frac{1}{i^{1-\alpha }}%
\sum\limits_{j=i}^{\left( n-1\right) /2}\frac{\lambda _{j}}{j^{1-\beta }}%
\frac{1}{\lambda _{j}}\right.
\end{equation*}
\begin{equation*}
+\left. \sum\limits_{j=\tau \left( n,m\right) }^{\left( n-1\right) /2}\frac{1%
}{j^{1-\beta }}\sum\limits_{i=j}^{\left( m-1\right) /2}\frac{1}{\lambda _{i}}%
\frac{\lambda _{i}}{i^{1-\alpha }}\right\}
\end{equation*}
\begin{equation*}
\left( \int\limits_{0}^{\pi /m}\int\limits_{0}^{\pi /n}|\Delta
_{ij}^{mn}f\left( x+s,y+t\right) |dsdt\right)
\end{equation*}
\begin{equation*}
\leq c\left( \alpha ,\beta \right) n\sum\limits_{i=\tau \left( n,m\right)
}^{\left( m-1\right) /2}\frac{\lambda _{i}}{i^{2-\left( \alpha +\beta
\right) }}\left( \int\limits_{0}^{\pi
/n}\sup\limits_{x}\sum\limits_{j=i}^{\left( n-1\right) /2}\frac{1}{\lambda
_{j}}|_{2}\Delta _{j}^{n}f\left( x,y+t\right) |dt\right)
\end{equation*}
\begin{equation*}
+c\left( \alpha ,\beta \right) m\sum\limits_{j=\tau \left( n,m\right)
}^{\left( n-1\right) /2}\frac{\lambda _{j}}{i^{2-\left( \alpha +\beta
\right) }}\left( \int\limits_{0}^{\pi
/m}\sup\limits_{y}\sum\limits_{i=j}^{\left( m-1\right) /2}\frac{1}{\lambda
_{j}}|_{1}\Delta _{i}^{m}f\left( x+s,y\right) |ds\right)
\end{equation*}
\begin{equation*}
\leq c\left( \alpha ,\beta \right) \left( V_{1}\Lambda \left( f\right)
+V_{2}\Lambda \left( f\right) \right) \sum\limits_{j=\tau \left( n,m\right)
}^{\infty }\frac{\lambda _{j}}{j^{2-\left( \alpha +\beta \right) }}%
\rightarrow 0,
\end{equation*}
as $n,m\rightarrow \infty$.

Combining (23)-(27) we conclude that
\begin{equation}
L_{mn}^{\left( 1\right) }\left( x,y\right) \rightarrow 0\,\,\,\,\,\text{%
as\thinspace \thinspace \thinspace }m,n\rightarrow \infty .
\end{equation}

From (18), (19), (21), (22) and (28) we obtain
\begin{equation}
L_{mn}\left( x,y\right) \rightarrow 0\,\,\,\,\,\text{as\thinspace \thinspace
\thinspace }m,n\rightarrow \infty .
\end{equation}

Finally, combining (7), (8), (16), (17) and (29) we get
\begin{equation*}
I_{mn}^{\left( 1\right) }\left( x,y\right) \rightarrow
0\,\,\,\,\,as\,\,\,m,n\rightarrow \infty .
\end{equation*}

Analogously, we can prove that
\begin{equation*}
I_{mn}^{\left( k\right) }\left( x,y\right) \rightarrow
0\,\,\,\,\,as\,\,\,m,n\rightarrow \infty ,\,\,k=2,3,4.
\end{equation*}
To complete the proof of Theorem 1, note that if $f$ is continuous on some compact $K$, then the relations $$
\lim_{s,t\to 0}\varphi_i(x,y,s,t)=0,\quad i=1,2,3,4,
$$
hold uniformly on $(x,y)\in K$ and all estimates in the proof also hold uniformly on $(x,y)\in K$. Hence
the $(C;-\alpha,\beta)$ means $\sigma _{n,m}^{\alpha ,\beta }(f;x,y)$ will converge to$f$ uniformly on $K$.
\end{proof}

\begin{proof}[Proof of Theorem 2]

It is not hard to see, that for any sequence $\Lambda =\{\lambda _{n}\}$
satisfying (1) the class $C\left( T^{2}\right) \bigcap P\Lambda BV$ is a
Banach space with the norm
\begin{equation*}
\left\| f\right\| _{P\Lambda BV}:=\left\| f\right\| _{C}+P\Lambda V\left(
f\right) .
\end{equation*}

Denote
\begin{equation*}
A_{i,j}:=\left[ \frac{\pi i-\alpha \pi /2}{N+1/2-\alpha /2},\frac{\pi \left(
i+1\right) -\alpha \pi /2}{N+1/2-\alpha /2}\right) \times \left[ \frac{\pi
j-\beta \pi /2}{N+1/2-\beta /2},\frac{\pi \left( j+1\right) -\beta \pi /2}{%
N+1/2-\beta /2}\right)
\end{equation*}
and
\begin{equation*}
W:=\left\{ \left( i,j\right) :j<i<2j,1<j<\frac{N-1}{2}\right\} .
\end{equation*}

Let
\begin{eqnarray*}
f_{N}\left( x,y\right) &:&=\sum\limits_{\left( i,j\right) \in W}t_{j}\mathbf{%
1}_{A_{i,j}}\left( x,y\right) \sin \left[ \left( N+1/2-\alpha /2\right)
x+\alpha \pi /2\right] \\
&&\times \sin \left[ \left( N+1/2-\beta /2\right) y+\beta \pi /2\right] ,
\end{eqnarray*}
where
\begin{equation*}
t_{j}:=\left( \sum\limits_{i=1}^{j}\frac{1}{\lambda _{i}}\right) ^{-1}.
\end{equation*}

First, we prove that $f\in P\Lambda BV.$ Indeed, let
\begin{equation*}
y\in \left[ \frac{\pi j-\beta \pi /2}{N+1/2-\beta /2},\frac{\pi \left(
j+1\right) -\beta \pi /2}{N+1/2-\beta /2}\right) .
\end{equation*}
Then it is evident that
\begin{equation*}
\sum\limits_{i}\frac{\left| f\left( I_{i},y\right) \right| }{\lambda _{i}}%
\leq c\left( \sum\limits_{i=j}^{2j-1}\frac{1}{\lambda _{2j-i}}\right)
t_{j}\leq c<\infty ,
\end{equation*}
consequently,
\begin{equation}
V_{1}\Lambda \left( f\right) <\infty .
\end{equation}

Let
\begin{equation*}
x\in \left[ \frac{\pi i-\alpha \pi /2}{N+1/2-\alpha /2},\frac{\pi \left(
i+1\right) -\alpha \pi /2}{N+1/2-\alpha /2}\right)
\end{equation*}
then from construction of the function $f$ we have
\begin{equation*}
\sum\limits_{j}\frac{\left| f\left( x,J_{j}\right) \right| }{\lambda _{j}}%
\leq c\sum\limits_{j=\left[ i/2\right] }^{i}\frac{t_{j}}{\lambda _{j-\left[
i/2\right] +1}}\leq ct_{\left[ i/2\right] }\left(
\sum\limits_{j=1}^{i-\left[ i/2\right] +1}\frac{1}{\lambda _{j}}\right) \leq
c<\infty .
\end{equation*}
Hence
\begin{equation}
V_{2}\Lambda \left( f\right) <\infty .
\end{equation}

Combining (30) and (31) and we conclude that $f\in P\Lambda BV.$

From (2)-(5) we can write
\begin{equation}
\pi ^{2}\sigma _{N,N}^{\left( -\alpha ,-\beta \right) }f_{N}\left( 0,0\right)
\end{equation}
\begin{equation*}
=\int\limits_{T^{2}}f_{N}\left( x,y\right) K_{N}^{-\alpha }\left( x\right)
K_{N}^{-\beta }\left( y\right) dxdy
\end{equation*}
\begin{equation*}
=\sum\limits_{\left( i,j\right) \in W}t_{j}\int\limits_{A_{i,j}}\sin \left[
\left( N+1/2-\alpha /2\right) x+\alpha \pi /2\right] \sin \left[ \left(
N+1/2-\beta /2\right) y+\beta \pi /2\right]
\end{equation*}
\begin{equation*}
\times O\left( \frac{1}{Nx^{2}}\right) O\left( \frac{1}{Ny^{2}}\right) dxdy
\end{equation*}
\begin{equation*}
+\sum\limits_{\left( i,j\right) \in W}t_{j}\int\limits_{A_{i,j}}\sin \left[
\left( N+1/2-\alpha /2\right) x+\alpha \pi /2\right] \frac{\sin ^{2}\left[
\left( N+1/2-\beta /2\right) y+\beta \pi /2\right] }{A_{N}^{-\beta }\left(
2\sin y/2\right)^{1-\beta} }
\end{equation*}
\begin{equation*}
\times O\left( \frac{1}{Nx^{2}}\right) dxdy
\end{equation*}
\begin{equation*}
+\sum\limits_{\left( i,j\right) \in W}t_{j}\int\limits_{A_{i,j}}\frac{\sin
^{2}\left[ \left( N+1/2-\alpha /2\right) x+\alpha \pi /2\right] }{%
A_{N}^{-\alpha }\left( 2\sin x/2\right)^{1-\alpha} }\sin \left[ \left( N+1/2-\beta
/2\right) y+\beta \pi /2\right]
\end{equation*}
\begin{equation*}
\times O\left( \frac{1}{Ny^{2}}\right) dxdy
\end{equation*}
\begin{equation*}
+\sum\limits_{\left( i,j\right) \in W}t_{j}\int\limits_{A_{i,j}}\frac{\sin
^{2}\left[ \left( N+1/2-\alpha /2\right) x+\alpha \pi /2\right] }{%
A_{N}^{-\alpha }\left( 2\sin x/2\right)^{1-\alpha} }\frac{\sin ^{2}\left[ \left(
N+1/2-\beta /2\right) y+\beta \pi /2\right] }{A_{N}^{-\beta }\left( 2\sin
y/2\right)^{1-\beta} }dxdy
\end{equation*}
\begin{equation*}
=:\sum\limits_{k=1}^{4}F_{N}^{\left( k\right) }\left( x,y\right)
\end{equation*}

It is easy to show that
\begin{equation}
\left| F_{N}^{\left( 1\right) }\left( x,y\right) \right| \leq
c\sum\limits_{\left( i,j\right) \in W}\frac{t_{j}}{ij}
\end{equation}
\begin{equation*}
=c\sum\limits_{j=1}^{\left[ \left( N-1\right) /2\right] }\frac{t_{j}}{j}%
\sum\limits_{i=j+1}^{2j-1}\frac{1}{i}\leq c\sum\limits_{j=1}^{\left[ \left(
N-1\right) /2\right] }\frac{t_{j}}{j},
\end{equation*}
\begin{equation}
\left| F_{N}^{\left( 2\right) }\left( x,y\right) \right| \leq c\left( \alpha
,\beta \right) \sum\limits_{\left( i,j\right) \in W}\frac{t_{j}}{ij^{1-\beta
}}
\end{equation}
\begin{equation*}
=c\left( \alpha ,\beta \right) \sum\limits_{j=1}^{\left[ \left( N-1\right)
/2\right] }\frac{t_{j}}{j^{1-\beta }}\sum\limits_{i=j+1}^{2j-1}\frac{1}{i}
\end{equation*}
\begin{equation*}
\leq c\left( \alpha ,\beta \right) \sum\limits_{j=1}^{\left[ \left(
N-1\right) /2\right] }\frac{t_{j}}{j^{1-\beta }},
\end{equation*}
\begin{equation}
\left| F_{N}^{\left( 3\right) }\left( x,y\right) \right| \leq c\left( \alpha
,\beta \right) \sum\limits_{j=1}^{\left[ \left( N-1\right) /2\right] }\frac{%
t_{j}}{j}\sum\limits_{i=j+1}^{2j-1}\frac{1}{i^{1-\alpha }}
\end{equation}
\begin{equation*}
\leq c\left( \alpha ,\beta \right) \sum\limits_{j=1}^{\left[ \left(
N-1\right) /2\right] }\frac{t_{j}}{j^{1-\alpha }}.
\end{equation*}

From the construction of the function $f_{N}$ we can write
\begin{equation}
\left| F_{N}^{\left( 4\right) }\left( x,y\right) \right|
\end{equation}
\begin{equation*}
=\frac{1}{\left( N+1/2-\alpha /2\right) \left( N+1/2-\beta /2\right) }%
\sum\limits_{\left( i,j\right) \in W}t_{j}
\end{equation*}
\begin{equation*}
\int\limits_{\pi i}^{\pi \left( i+1\right) }\int\limits_{\pi j}^{\pi \left(
j+1\right) }\frac{\sin ^{2}u}{A_{N}^{-\alpha }\left( 2\sin \frac{u-\alpha
\pi /2}{2\left( N+1/2-\alpha /2\right) }\right) ^{1-\alpha }}\frac{\sin ^{2}v%
}{A_{N}^{-\alpha }\left( 2\sin \frac{v-\beta \pi /2}{2\left( N+1/2-\beta
/2\right) }\right) ^{1-\beta }}dudv
\end{equation*}
\begin{equation*}
\geq \frac{c\left( \alpha ,\beta \right) N^{\alpha +\beta }}{N^{2}}%
\sum\limits_{\left( i,j\right) \in W}t_{j}\frac{N^{2-\left( \alpha +\beta
\right) }}{i^{1-\alpha }j^{1-\beta }}\int\limits_{\pi i}^{\pi \left(
i+1\right) }\sin ^{2}udu\int\limits_{\pi j}^{\pi \left( j+1\right) }\sin
^{2}vdv
\end{equation*}
\begin{equation*}
\geq c\left( \alpha ,\beta \right) \sum\limits_{\left( i,j\right) \in W}%
\frac{t_{j}}{j^{1-\beta }}\frac{1}{i^{1-\alpha }}
\geq c\left( \alpha ,\beta \right) \sum\limits_{j=1}^{\left[ \left(
N-1\right) /2\right] }\frac{t_{j}}{j^{1-\beta }}\sum\limits_{i=j+1}^{2j-1}%
\frac{1}{i^{1-\alpha }}
\end{equation*}
\begin{equation*}
\geq c\left( \alpha ,\beta \right) \sum\limits_{j=1}^{\left[ \left(
N-1\right) /2\right] }\frac{t_{j}}{j^{1-\left( \beta +\alpha \right) }}.
\end{equation*}

Since $\frac 1 {j^{1-\alpha}}+\frac 1 {j^{1-\beta}}=o\left(\frac 1 {j^{1-(\alpha+\beta)}}\right)$ as $j\to\infty$,
from (32)-(36) we conclude that if $j_0$ is big enough and $N>2j_0$, then
\begin{equation}
\pi ^{2}\left| \sigma _{N,N}^{\left( -\alpha ,-\beta \right) }f_{N}\left(
0,0\right) \right| \geq c\left( \alpha ,\beta \right)
\sum\limits_{j=j_{0}}^{\left[ \left( N-1\right) /2\right] }\frac{t_{j}}{%
j^{1-\left( \beta +\alpha \right) }}.
\end{equation}

Let $\lambda _{j}=j^{1-\left( \alpha +\beta \right) }\gamma
_{j},\,\,\,\gamma _{j}\geq \gamma _{j+1},\,\,j=1,2,....$ Then we can write
\begin{equation*}
\frac{1}{t_{j}}=\sum\limits_{i=1}^{j}\frac{1}{\lambda _{i}}%
=\sum\limits_{i=1}^{j}\frac{1}{i^{1-\left( \alpha +\beta \right) }\gamma _{i}%
}\leq c\left( \alpha ,\beta \right) \frac{j^{\alpha +\beta }}{\gamma _{j}}.
\end{equation*}
Consequently,
\begin{equation}
t_{j}j^{\alpha +\beta }\geq c\left( \alpha ,\beta \right) \gamma _{j}.
\end{equation}

Combining (37) and (38) we obtain
\begin{equation}
\pi ^{2}\left| \sigma _{N,N}^{\left( -\alpha ,-\beta \right) }f_{N}\left(
0,0\right) \right| \geq c\left( \alpha ,\beta \right)
\sum\limits_{j=j_{0}}^{\left[ \left( N-1\right) /2\right] }\frac{\gamma _{j}%
}{j}
\end{equation}
\begin{equation*}
=c\left( \alpha ,\beta \right) \sum\limits_{j=j_{0}}^{\left[ \left(
N-1\right) /2\right] }\frac{t_{j}}{j^{2-\left( \beta +\alpha \right) }}%
\rightarrow \infty \text{\thinspace \thinspace \thinspace as\thinspace
\thinspace \thinspace }N\rightarrow \infty .
\end{equation*}

Applying the Banach-Steinhaus Theorem, from (39) we obtain that there exists
a continuous function $f\in P\Lambda BV$ such that
\begin{equation*}
\sup\limits_{N}\left| \sigma _{N,N}^{\left( -\alpha ,-\beta \right)
}f_{N}\left( 0,0\right) \right| =+\infty .
\end{equation*}
\end{proof}
{\bf Acknowledgement}. The authors would like to thank the referee for helpful suggestions.

\end{document}